\documentclass[12pt]{amsart}
\usepackage{amssymb}
\usepackage{bbm}
\usepackage{amsfonts}
\usepackage{latexsym}
\usepackage[all]{xy}
\usepackage{amscd}
\usepackage{comment}
\usepackage{fullpage}

\newtheorem{theorem}{Theorem}[section]
\newtheorem{lemma}[theorem]{Lemma}

\newtheorem{corollary}[theorem]{Corollary}
\theoremstyle{definition}

\newtheorem{example}[theorem]{Example}

\newtheorem{remark}[theorem]{Remark}

\newcommand{\DOT}{\setlength{\unitlength}{1pt}\begin{picture}(2.5,2)
                  (1,1)\put(2,3.5){\circle*{2}}\end{picture}}

\newcommand{\bu}{\DOT}
\newcommand{\la}{\langle}
\newcommand{\ra}{\rangle}
\newcommand{\ot}{\otimes}
\newcommand{\Wedge}{\textstyle\bigwedge}
\newcommand{\lexp}[2]{{\vphantom{#2}}^{#1}{#2}}
\renewcommand{\k}{\Bbbk}
\newcommand{\N}{\mathbb{N}}
\newcommand{\M}{\mathbb{N}\cup\{-1\}}
\renewcommand{\O}{\{0,1\}}
\newcommand{\HHD}{\HH^{\bu}}
\newcommand{\Ph}{\Theta}

\DeclareMathOperator{\Hom}{Hom}
\DeclareMathOperator{\op}{op} 
\DeclareMathOperator{\HH}{HH}
\DeclareMathOperator{\Ext}{Ext}
\DeclareMathOperator{\sgn}{sgn}

\renewcommand{\span}{\operatorname{span}}

\begin{document}

\title[Hochschild cohomology]
{Hochschild cohomology of group extensions of quantum symmetric algebras}
\date{May 14, 2010}


\author{Deepak Naidu}
\email{dnaidu@math.tamu.edu}
\author{Piyush Shroff}
\email{pshroff@math.tamu.edu}
\author{Sarah Witherspoon}
\email{sjw@math.tamu.edu}
\address{Department of Mathematics, Texas A\&M University,
College Station, Texas 77843, USA}
\thanks{The second and third authors were partially supported by
NSF grant \#DMS-0800832 and Advanced Research Program Grant 010366-0046-2007
from the Texas Higher Education Coordinating Board.}

\begin{abstract}
Quantum symmetric algebras (or noncommutative polynomial rings)
arise in many places in mathematics. In this article we find the
multiplicative structure of their Hochschild cohomology when the
coefficients are in an arbitrary bimodule algebra. When this bimodule
algebra is a finite group extension (under a diagonal action) of a quantum
symmetric algebra, we give explicitly the graded vector space structure.
This yields a complete description of the Hochschild cohomology ring of
the corresponding skew group algebra.
\end{abstract}

\maketitle

\begin{section}{Introduction}
The action of a group $G$ on a vector space $V$ induces an action
on the corresponding symmetric algebra $S(V)$ (a polynomial ring).
The resulting skew group algebra $S(V)\rtimes G$ is a noncommutative
ring encoding the group action. 
Deformations of such skew group algebras are related to deformations of
corresponding orbifolds, and some have appeared independently in several
places under different names: graded Hecke algebras, rational Cherednik
algebras, and symplectic reflection algebras. Deformations of any algebra are
intimately related to its Hochschild cohomology. When $G$ is finite
the Hochschild cohomology of $S(V)\rtimes G$ was
computed independently by Farinati \cite{F} and by Ginzburg and Kaledin
\cite{GK}. Its algebra structure was first given by Anno \cite{A}.

In this paper we replace the symmetric algebra $S(V)$ with a quantum symmetric
algebra and explore its Hochschild cohomology.
Our quantum symmetric algebra is a noncommutative
polynomial ring, denoted $S_{\bf q}(V)$, 
in which the variables commute only up to
multiplication by nonzero scalars (encoded in the vector ${\bf q}$).
Noncommutative polynomials have been of interest for some time.
Our current work incorporates group actions as well and is in part
motivated by the recent appearance of two articles: 
Kirkman, Kuzmanovich, and Zhang \cite{KKZ} prove a version of the 
classical Shephard-Todd-Chevalley Theorem, namely that the invariant
subring of $S_{\bf q}(V)$ under a finite group action is again a 
quantum symmetric algebra. In the setting of an ordinary polynomial
ring, methods from invariant theory for finding such subrings
play a crucial role in computations of Hochschild cohomology (see,
for example, \cite{SW}). 
Around the same time, Bazlov and Berenstein \cite{BB} introduced
braided Cherednik algebras, which are deformations of $S_{\bf q}(V)\rtimes G$
in some special cases.
Knowledge of Hochschild cohomology will provide insight into
these and other possible deformations.

More specifically, let $\k$ be a field of characteristic 0.
Let $N$ be a positive integer and for each pair $i,j$
of elements in $\{1,\ldots,N\}$, let $q_{i,j}$ be a nonzero scalar such that
$q_{i,i}=1$ and $q_{j,i}=q_{i,j}^{-1}$ for all $i,j$. Denote by ${\bf q}$
the corresponding tuple of scalars, 
${\bf q} := (q_{i,j})_{1\leq i < j \leq N}$.
Let $V$ be a vector space with basis $x_1,\ldots,x_N$, and 
let 
$$
 S_{\bf q}(V) := \k\la x_1,\ldots,x_N \mid x_ix_j = q_{i,j}x_jx_i  \mbox{ for all }
      1\leq i,j\leq N \ra,
$$
the {\em quantum symmetric algebra} determined by $\bf q$.
This is a Koszul algebra (see e.g.\ \cite{AS}), so there is a standard complex
$K_{\bu}(S_{\bf q}(V))$ that is a resolution of $S_{\bf q}(V)$ as an
$S_{\bf q}(V)$-bimodule. This complex is given in Section \ref{free-resolution}
below; see \cite{W} for details for more general quantum symmetric
algebras arising from braidings. 
Priddy \cite{P} first introduced Koszul algebras, and the theory was
developed further, including such complexes, in \cite{BGS,BG,M}.

We first compute cup products, using the resolution in Section \ref{free-resolution}, 
and obtain the following theorem.
Here $\Wedge_{{\bf q}^{-1}}(V^*)$ denotes a quantum exterior algebra,
defined in (\ref{qea}) below, on the dual vector space $V^*$, and 
$B$ is any algebra with a compatible structure of an
$S_{\bf q}(V)$-bimodule. The two choices of algebra $B$ in which we
are most interested are $B=S_{\bf q}(V)$ and $B=S_{\bf q}(V)\rtimes G$,
where $G$ is a finite group of graded automorphisms of $S_{\bf q}(V)$.
 
\medskip

\noindent{\bf Theorem \ref{cup product HH(A,B)}.} \textit{The
Hochschild cohomology $\HHD(S_{\bf q}(V), B)$
is a subquotient algebra of the tensor product 
$B\ot \Wedge_{{\bf q}^{-1}}(V^*)$.}

\medskip

That is, the Hochschild cohomology is a vector subquotient of $B
\ot \Wedge_{{\bf q}^{-1}}(V^*)$,
and its cup product is determined by that of the tensor product of
the two algebras $B$ and $\Wedge_{{\bf q}^{-1}}(V^*)$. 

We next give the graded vector space structure of Hochschild cohomology 
in the two cases $B=S_{\bf q}(V)$ and $B=S_{\bf q}(V)\rtimes G$, when $G$
acts diagonally on the basis $x_1,\ldots,x_N$ of $V$.
We adapt techniques developed in the context of Hochschild homology by Wambst \cite{W}.
In the first case this gives the following technical result.
The notation will be explained in Section~\ref{gvsg}.

\medskip

\noindent
{\bf Corollary \ref{HH(A)}.} \textit{For all $m\in \N$,
$$
\HH^m(S_{\bf q}(V)) \cong \bigoplus_{\substack{\beta \in \O^N \\ |\beta| = m}} 
\bigoplus_{\substack{\alpha \in \N^N \\ \alpha - \beta \in C}}
\span_\k\{x^\alpha \ot {(x^*)}^{\wedge \beta}\}.
$$
}

\medskip

This result is a consequence of our more general Theorem \ref{HH(A:G)} on $S_{\bf q}(V)\rtimes G$.
It should be compared with work of Richard \cite{R}, in which there are
some results on the Hochschild cohomology 
of a related ring of twisted differential operators on quantum affine space.
Richard obtains his results by first computing Hochschild homology and then 
invoking a duality between homology and cohomology.
In our setting, such a duality does not hold in general; a comparison of our
Example~\ref{N=2} below with Wambst's Corollary 6.2 \cite{W} shows that there
is no duality in our smallest possible case. Other authors have used various
techniques to compute Hochschild homology of generalizations of
$S_{\bf q}(V)$ \cite{GG,Gu,W}.

Since the characteristic of $\k$ is 0, the Hochschild cohomology of $S_{\bf q}(V)\rtimes G$
is the subalgebra of $G$-invariant elements of $\HH^{\bu}(S_{\bf q}(V), S_{\bf q}(V)\rtimes G)$,
and the multiplicative structure of this algebra is given by Theorem \ref{cup product HH(A,B)}.
To determine the additive structure precisely, we specialize to diagonal actions of $G$
on the chosen basis $x_1,\ldots,x_N$ of $V$.
This is the case to which Wambst's techniques in \cite{W} may be most easily
adapted to express the relevant cochain complex as the direct sum of an acyclic
complex and a complex in which all differentials are~0.
This leads to our description of the Hochschild cohomology as a graded
vector space in the following theorem.
The notation will be explained in Section~\ref{gvsg}.

\medskip

\noindent
{\bf Theorem \ref{HH(A:G)}.} \textit{Assume the finite group $G$ acts diagonally on the chosen
basis of $V$. Then
for all $m\in \N$, 
$$
\HH^m(S_{\bf q}(V),S_{\bf q}(V) \rtimes G) 
\cong \bigoplus_{g \in G}
\bigoplus_{\substack{\beta \in \O^N \\ |\beta| = m}} 
\bigoplus_{\substack{\alpha \in \N^N \\ \alpha - \beta \in C_g}}
\span_\k\{(x^\alpha\#g) \ot {(x^*)}^{\wedge \beta}\},
$$
and $\HH^m(S_{\bf q}(V)\rtimes G)$ is its $G$-invariant subspace. 
}

\medskip

We plan to address more general group actions, as well as related
deformations of the skew group algebra $S_{\bf q}(V)\rtimes G$, in future articles.\\

\noindent {\bf Organization.}
This paper is organized as follows. Section~\ref{prelim} contains necessary preliminary
information, including the  resolution of $S_{\bf q}(V)$ that will be used.
We also give a chain map from this resolution to the bar
resolution of $S_{\bf q}(V)$, used in Section~\ref{B} to compute cup products.

We prove Theorem \ref{cup product HH(A,B)} in Section~\ref{B}.
In Section~\ref{gvsg} we prove Theorem~\ref{HH(A:G)} and apply 
Theorem~\ref{cup product HH(A,B)} to give the 
cup product on $\HHD(S_{\bf q}(V)\rtimes G)$.
As a special case, when $G=\{1\}$, we obtain in Corollary~\ref{HH(A)} the
Hochschild cohomology of $S_{\bf q}(V)$.

\end{section}

\begin{section}{Preliminaries}
\label{prelim}

All tensor products and exterior powers are taken over the field $\k$
of characteristic 0.

Let $A$ be an algebra over $\k$, and let $M$ be an $A$-bimodule.
We identify $M$ with a (left) $A^e$-module, where $A^e=A\ot A^{\op}$;
$A^{\op}$ is the algebra $A$ with the opposite multiplication.
The {\em Hochschild cohomology} of $A$ with coefficients in $M$ is
$$
  \HH^{\bu}(A,M) := \Ext^{\bu}_{A^e}(A,M),
$$
where $A$ is considered to be an $A^e$-module under left and right multiplication.
One useful free $A^e$-resolution of $A$ is the {\em bar resolution}
\begin{equation}\label{bar}
 \cdots \xrightarrow{\delta_3} A^{\ot 4} \xrightarrow{\delta_2}
  A^{\ot 3}\xrightarrow{\delta_1} A^e \xrightarrow{\text{mult}} A \xrightarrow{ } 0 
\end{equation}
where 
$\delta_m(a_0\ot \cdots\ot a_{m+1}) =\sum_{i=0}^m (-1)^i a_0\ot\cdots\ot
    a_ia_{i+1}\ot \cdots \ot a_{m+1}$ for all $a_0,\ldots,a_{m+1}\in A$, and
the map from $A^e$ to $A$ is given by multiplication in $A$.
Suppose $M=B$ is an $A$-{\em bimodule algebra}, 
that is $B$ is an algebra and also an $A$-bimodule and these two structures are compatible
in the sense that $a(bb') = (ab)b'$  and $(bb')a = b (b'a)$ for all $a\in A$ and $b,b'\in B$.
Then $\HH^{\bu}(A,B)$ has a cup product defined at the cochain level as follows (e.g.\
see \cite[p.\ 278]{G}).
Let $f\in \Hom_{A^e}(A^{\ot (m+2)}, B)$, $f'\in \Hom_{A^e}(A^{(n+2)},B)$.
Then $f\smile f'\in \Hom_{A^e}(A^{\ot (m+n+2)}, B)$ is determined by
$$
  f\smile f' (a_0\ot a_1\ot\cdots \ot a_{m+n}\ot a_{m+n+1}) =
 f(a_0\ot a_1\ot\cdots\ot a_m\ot 1) f'(1\ot a_{m+1}\ot \cdots\ot a_{m+n}\ot a_{m+n+1}).
$$

Let $G$ be a finite group acting on the algebra $A$ by automorphisms.
We denote by ${}^ga$ the result of applying $g\in G$ to $a\in A$.
Then we may form the {\em skew group algebra} $A\rtimes G$:
Additively, it is the free $A$-module with basis $G$.
We write $A\rtimes G = \oplus_{g\in G}A_g$,
where $A_g = \{a\# g\mid a\in A\}$, that is for each $a\in A$ and
$g\in G$ we denote by $a\# g\in A_g$ the $a$-multiple of~$g$.
Multiplication on $A\rtimes G$ is determined by 
$$
   (a\# g) (b\# h) := a ({}^gb) \# gh
$$
for all $a,b\in A$, $g,h\in G$. Note that for each $g\in G$, $A_g$ is a (left)
$A^e$-module via the action 
$$
   (a\ot b) \cdot (c\# g) := (a\# 1)(c\# g) (b\# 1) = ac ({}^g b) \# g
$$
for all $a,b,c\in A$, $g\in G$.

For convenience in what follows, we will sometimes denote the quantum symmetric algebra 
$S_{\bf q}(V)$ simply by $A$.

\begin{subsection}{A free resolution of $S_{\bf q}(V)$}
\label{free-resolution}

By \cite[Proposition 4.1(c)]{W},
the following is a free $A^e$-resolution of $A=S_{\bf q}(V)$:

\begin{equation}
\label{label: free resolution}
\cdots \xrightarrow{} A^e\ot\Wedge^2(V) \xrightarrow{d_2}
A^e \otimes \Wedge ^1(V) \xrightarrow{d_1}
A^e \xrightarrow{\text{mult}} A \xrightarrow{} 0,
\end{equation}
that is, for $1\leq m\leq N$, the degree $m$ term is $A^e\otimes \Wedge^m(V)$, and 
 $d_m$  is defined by 
\begin{equation*}
\begin{split}
&d_m(1^{\ot 2}\ot x_{j_1}\wedge\cdots\wedge x_{j_m})\\
&=  \sum_{i=1}^m (-1)^{i+1} \left[ \left(\prod_{s=1}^{i} q_{j_s, j_i} \right)
    x_{j_i}\ot 1 - \left(\prod_{s=i}^m q_{j_i, j_s} \right) \ot x_{j_i} \right] \ot
    x_{j_1}\wedge \cdots \wedge \hat{x}_{j_i} \wedge\cdots \wedge x_{j_m}
\end{split}
\end{equation*}
whenever $1\leq j_1 <\ldots < j_m\leq N$.
This is a twisted version of the usual Koszul resolution for a polynomial ring.

Let us write the above formula for $d_m$ in a more convenient form. 
We first introduce some notation following  Wambst \cite{W}. 
Let $\N^N$ denote the set of all $N$-tuples of elements
from $\N$. For any $\alpha \in \N^N$, the {\em length}
of $\alpha$, denoted $|\alpha|$, is the sum $\sum_{i=1}^N \alpha_i$.
For all $\alpha \in \N^N$, define $x^\alpha := 
x_1^{\alpha_1} x_2^{\alpha_2} \cdots x_N^{\alpha_N}$. For all 
$i \in \{1, \ldots, N\}$,  define $[i] \in \N^N$ by
$[i]_j = \delta_{i,j}$, for all $j \in \{1, \ldots, N\}$.
For any $\beta \in \O^N$, let $x^{\wedge \beta}$ denote
the vector $x_{j_1} \wedge \cdots \wedge x_{j_m} \in \Wedge^m(V)$
which is defined by $m = |\beta|$, $\beta_{j_k} = 1$ for all 
$k \in \{1, \ldots, m\}$, and $j_1<\ldots <j_m$.
Then, for any $\beta \in \O^N$ with $|\beta| = m$ we have
\begin{equation*}
d_m(1^{\ot 2}\ot x^{\wedge \beta}) =
\sum_{i=1}^N \delta_{\beta_i,1} (-1)^{\sum_{s=1}^{i-1} \beta_s} 
\left[ \left(\prod_{s=1}^i q_{s, i}^{\beta_s} \right)
    x_i\ot 1 - \left(\prod_{s=i}^N q_{i, s}^{\beta_s} \right) \ot x_i \right] \ot
    x^{\wedge(\beta-[i])}.
\end{equation*}

\end{subsection}

\begin{subsection}{A chain map into the bar resolution of $S_{\bf q}(V)$}

We wish to define a chain map from our complex $A^e\ot \Wedge^{\bu}(V)$ to the bar
resolution (\ref{bar}) for $A=S_{\bf q}(V)$.
Wambst defined a more general chain map \cite[Lemma 5.3 and Theorem 5.4]{W}.
Here we introduce notation useful in our setting, and include some details for completeness.

For each set of $m$ distinct natural numbers $j_1,\ldots, j_m$ ($m\leq N$)
and each permutation $\pi\in S_m$, the scalar $q_{\pi}^{j_1,\ldots,j_m}$ is
defined by the equation
$$
   q_{\pi}^{j_1,\ldots,j_m} x_{j_{\pi(1)}}\cdots x_{j_{\pi(m)}} 
  = x_{j_1}\cdots x_{j_m}
$$
in $A$, that is, $q_{\pi}^{j_1,\ldots,j_m}$ is the scalar arising when one applies
the permutation $\pi$ to the variables in the product $x_{j_1}\cdots x_{j_m}$,
using the relations in $A$ to rewrite it.  

The following lemma is immediate from the definition.

\begin{lemma}\label{sigma-tau}
If $\pi = \sigma\tau$ in $S_m$, then
$$
   q_{\pi}^{j_1,\ldots,j_m} = q_{\sigma}^{j_{\tau(1)},\ldots,j_{\tau(m)}} 
    q_{\tau}^{j_1,\ldots,j_m}.
$$
\end{lemma}

For each $m\geq 1$, define the map $\phi_m: A^e\ot \Wedge^m(V)\rightarrow A^{\ot (m+2)}$
by 
\begin{equation}\label{phim}
  \phi_m(1^{\ot 2}\ot x_{j_1}\wedge \cdots \wedge x_{j_m})
  = \sum_{\pi\in S_m} (\sgn\pi) q_{\pi}^{j_1,\ldots,j_m} \ot x_{j_{\pi(1)}}\ot \cdots
    \ot x_{j_{\pi(m)}} \ot 1
\end{equation}
for all distinct $x_{j_1},\ldots,x_{j_m}$.
Note that $\phi_m$ is injective: The image of a basis of $\Wedge^m(V)$ under $\phi_m$
is linearly independent, as may be seen by comparing the variables involved.
Set $\phi_0$ and $\phi_{-1}$ to be the identity maps on $A\ot A$ and $A$, respectively.

\begin{remark}\label{BG}
By its definition, the image of $\phi_m$ is contained in
\begin{equation}\label{R}
   \bigcap_{i=0}^{m-2} (A\ot V^{\ot i}\ot R\ot V^{\ot (m-i-2)} \ot A),
\end{equation}
where $R\subset V\ot V$ is the vector subspace spanned 
by the relations $x_i\ot x_j - q_{i,j}x_j\ot x_i$.
For example, to see that $\phi_m(1^{\ot 2}\ot x_{j_1}\wedge \cdots\wedge x_{j_m})$
is in $A\ot R\ot V^{\ot (m-2)}\ot A$, fix $\pi\in S_m$. Let $(12)$ denote
the permutation transposing 1 and 2.
The $\pi$- and $\pi (12)$-terms of formula (\ref{phim}) above are
(writing $\pi(12)= (\pi(1),\pi(2))\pi$ in $S_m$ and applying the lemma):
\begin{eqnarray*}
 &&\hspace{-1cm} (\sgn\pi)q_{\pi}^{j_1,\ldots,j_m} \ot x_{j_{\pi(1)}}\ot\cdots
   \ot x_{j_{\pi(m)}}\ot 1 + (\sgn \pi(12)) q_{\pi (12)}^{j_1,\ldots,j_m}
    \ot x_{j_{\pi(2)}}\ot x_{j_{\pi(1)}}\ot \cdots\ot x_{j_{\pi(m)}}\ot 1\\
 &=& (\sgn\pi) q_{\pi}^{j_1,\ldots,j_m}\ot x_{j_{\pi(1)}}\ot\cdots\ot
    x_{j_{\pi(m)}}\ot 1 \\
    &&- (\sgn\pi)
     q_{(\pi(1), \pi(2))}^{j_{\pi(1)},\ldots, j_{\pi(m)}} 
   q_{\pi}^{j_1,\ldots,j_m} 
   \ot x_{j_{\pi(2)}}
   \ot x_{j_{\pi(1)}}\ot \cdots\ot x_{j_{\pi(m)}}\ot 1\\
 &=& (\sgn\pi) q_{\pi}^{j_1,\ldots,j_m}\ot (x_{j_{\pi(1)}}\ot x_{j_{\pi(2)}}
    - q_{(\pi(1),\pi(2))}^{j_{\pi(1)},\ldots,j_{\pi(m)}} x_{j_{\pi(2)}}
   \ot x_{j_{\pi(1)}} ) \ot x_{j_{\pi(3)}}\ot\cdots\ot x_{j_{\pi(m)}}\ot 1,
\end{eqnarray*}
which is visibly in $A\ot R\ot V^{\ot (m-2)}\ot A$.
\end{remark}

\begin{lemma}
The map $\phi$ defined in equation \eqref{phim} is a chain map.
\end{lemma}

\begin{proof}
We must show that 
\begin{equation}\label{commutes}
\phi_{m-1}\circ d_m (1^{\ot 2}\ot x_{j_1}\wedge\cdots \wedge x_{j_m})
=\delta_m\circ\phi_m(1^{\ot 2}\ot x_{j_1}\wedge\cdots \wedge x_{j_m})
\end{equation}
for all $j_1,\ldots, j_m$ with $j_1<\cdots < j_m$, 
where $\delta_m$ is the differential on  
the bar complex (\ref{bar}).
This  may easily be checked when $m=0$.

By Remark \ref{BG}, the formula (\ref{phim}) for $\phi_m$ and the formula for $\delta_m$,
the right side of (\ref{commutes}) is
$$
    \sum_{\pi\in S_m} (\sgn\pi) q_{\pi}^{j_1,\ldots,j_m} \left(x_{j_{\pi(1)}}\ot \cdots 
    \ot x_{j_{\pi(m)}}\ot 1 
   + (-1)^m \ot 
     x_{j_{\pi(1)}}\ot \cdots \ot x_{j_{\pi(m)}}\right)
$$
for all $m\geq 1$.
The left side of (\ref{commutes}) is 
$$
\begin{aligned}
  & \phi_{m-1}\left(\sum_{i=1}^m (-1)^{i+1} \left( (\prod_{s=1}^i q_{j_s,j_i})x_{j_i}\ot 1 -
    (\prod_{s=i}^m q_{j_i,j_s})\ot x_{j_i}\right) \ot x_{j_1}\wedge\cdots\wedge
   \hat{x}_{j_i}\wedge\cdots\wedge x_{j_m} \right)\\
 &= \sum_{i=1}^m (-1)^{i+1}  \left( (\prod_{s=1}^i q_{j_s,j_i})x_{j_i}\ot 1 -
    (\prod_{s=i}^m q_{j_i,j_s})\ot x_{j_i}\right)
   \\
   & \hspace{3cm} \ot \left(\sum_{\pi '\in S^i_{m-1}} (\sgn\pi ') q_{\pi'}^{j_1,\ldots,
    \hat{j_i},
    \ldots,j_m} x_{j_{\pi'(1)}}\ot\cdots\ot x_{j_{\pi'(m)}}\right) \\
  & = \sum_{i=1}^m (-1)^{i+1} \sum_{\pi'\in S^i_{m-1}} 
   \left( (\sgn\pi')(\prod_{s=1}^i q_{j_s,j_i}) q_{\pi'}^{j_1,\ldots,\hat{j_i},\ldots,j_m}
    x_{j_i}\ot x_{j_{\pi'(1)}}\ot\cdots\ot x_{j_{\pi'(m)}}\ot 1 \right.\\
   & \hspace{3cm} - \left.  (\sgn\pi')(\prod_{s=i}^m q_{j_i,j_s})
   q_{\pi'}^{j_1,\ldots,\hat{j_i},\ldots,j_m} \ot x_{j_{\pi'(1)}}\ot \cdots\ot
    x_{j_{\pi'(m)}}\ot x_{j_i} \right),
\end{aligned}
$$
where for each $i$, $S_{m-1}^i$ denotes the symmetric group on $\{1,\ldots, \hat{i},\ldots, m\}$.
In the first set of summands  above, if
$
  \pi := \pi' (i, \ i-1, \ \ldots, \ 1),
$
then $\sgn\pi ' = (-1)^{i+1} \sgn\pi$, and 
we may replace $(\sgn\pi')x_{j_i}\ot x_{j_{\pi'(1)}}\ot\cdots \ot x_{j_{\pi'(m)}}\ot 1$ by
$ (-1)^{i-1}(\sgn\pi) x_{j_{\pi(1)}}\ot \cdots\ot x_{j_{\pi(m)}}\ot 1$.
Similarly in the second set of summands,  let $\pi = \pi' (i, \ i+1, \ \ldots, \ m)$, 
and replace
$(\sgn\pi') \ot x_{j_{\pi'(1)}}\ot \cdots \ot x_{j_{\pi'(m)}}\ot x_{j_i}$
by $(-1)^{m-i} (\sgn\pi) \ot x_{j_{\pi(1)}}\ot\cdots\ot x_{j_{\pi(m)}}$.
Again, for the first set of summands, notice that
$q_{\pi}^{j_1,\ldots,j_m} = (\prod_{s=1}^i q_{j_s,j_i}) q_{\pi'}^{j_1,\ldots,\hat{j_i},
   \ldots,j_m}$, and for the second set,
$q_{\pi}^{j_1,\ldots,j_m} = (\prod_{s=i}^m q_{j_i,j_s}) q_{\pi'}^{j_1,\ldots,
  \hat{j_i},\ldots,j_m}$.
Making all such replacements, we find that the left side of
(\ref{commutes}) is indeed equal to the right side.
\end{proof}

\begin{remark}\label{complement}
The image of $\phi_m$ is a free $A^e$-submodule of $A^{\ot (m+2)}$ that is a direct
summand of $A^{\ot (m+2)}$ as an $A^e$-module: 
Take a vector space complement in $\k \ot V^{\ot m}\ot \k$ to the image of
$\k \ot \k \ot \Wedge^m(V)$ under $\phi_m$, and extend to the required
complementary $A^e$-module direct summand of $A^{\ot (m+2)}$.
It follows that there is a chain map $\psi$ from the bar resolution to 
$A^e\ot \Wedge^{\DOT} (V)$ for which $\psi\phi$ is the identity map. 
\end{remark}

\end{subsection}

\end{section}

\begin{section}{The cup product}
\label{B}

As before, let $A$ denote the quantum symmetric algebra $S_{\bf q}(V)$
and let $A^e$ denote the enveloping algebra $A \ot A^{\op}$. Let $B$ denote
an $A$-bimodule. In this section, we will describe the cup product on Hochschild cohomology 
$\HH^{\bu}(A,B)=\Ext^{\bu}_{A^e}(A,B)$, when $B$ additionally has a compatible algebra structure.

We begin by applying $\Hom_{A^e}(\cdot, B)$ to \eqref{label: free resolution},
dropping the term $\Hom_{A^e}(A,B)$ to obtain the complex
\begin{equation}
\label{Hom(resolution) for B}
0 \xrightarrow{} \Hom_{A^e}(A^e, B) \xrightarrow{\tilde{d}_1} 
\Hom_{A^e}(A^e \ot \Wedge^1(V), B) \xrightarrow{\tilde{d}_2} 
\Hom_{A^e}(A^e \ot \Wedge^2(V), B) \xrightarrow{} \cdots,
\end{equation}
that is, the degree $m$ term is $\Hom_{A^e}(A^e\otimes \Wedge^{m}(V),B)$, and 
 $\tilde{d}_m$ is defined by $\tilde{d}_m(F) = F \circ d_m$, for all
$F \in \Hom_{A^e}(A^e\otimes \Wedge^{m-1}(V),B)$.  

We identify the spaces $\Hom_{A^e}(A^e \ot \Wedge^m(V), B)$ and $B \ot \Wedge^m(V^*)$
via the map
$$ 
\Ph_m: \Hom_{A^e}(A^e \ot \Wedge^m(V), B) \to B \ot \Wedge^m(V^*)
$$
$$
\hspace{.8in} F \mapsto \displaystyle \sum_{\substack{\beta \in \O^N \\ |\beta|=m}}
F(1^{\ot 2} \ot x^{\wedge \beta}) \ot
{(x^*)}^{\wedge \beta}.
$$

We thus obtain the following complex,
which is equivalent to the complex in \eqref{Hom(resolution) for B}:
\begin{equation}
\label{new Hom(resolution) for B}
0 \xrightarrow{}  B \xrightarrow{d_1^*} 
B \ot \Wedge^1(V^*) \xrightarrow{d_2^*} 
B \ot \Wedge^2(V^*) \xrightarrow{} \cdots,
\end{equation}
that is, the degree $m$ term is $B \otimes \Wedge^{m}(V^*)$, and 
 $d_m^*$ is defined by $d_m^* = \Ph_m \circ \tilde{d}_m \circ \Ph^{-1}_{m-1}$,
where $\Ph^{-1}_0$ takes $b\in B$ to $F\in\Hom_{A^e}(A^e,B)$ defined by
$F(1^{\ot 2})=b$.

Let us describe $d_m^*$ explicitly. For any $b \in B$ and  
$\beta \in \O^N$ with $|\beta|=m-1$ we have
\begin{equation*}
\begin{split}
d_m^*&(b \ot {(x^*)}^{\wedge \beta})\\
&= (\Ph_m \circ \tilde{d}_m \circ \Ph^{-1}_{m-1})
(b \ot {(x^*)}^{\wedge \beta}) \\
&= \Ph_m(\Ph^{-1}_{m-1}(b \ot {(x^*)}^{\wedge \beta}) \circ d_m)\\
&= \sum_{\substack{\beta' \in \O^N \\ |\beta'|=m}}
(\Ph^{-1}_{m-1}(b \ot {(x^*)}^{\wedge \beta}) \circ d_m)(1^{\ot 2} \ot x^{\wedge \beta'}) 
\ot {(x^*)}^{\wedge \beta'}\\
&= \sum_{i=1}^N \delta_{\beta_i,0}
(\Ph^{-1}_{m-1}(b \ot {(x^*)}^{\wedge \beta}) \circ d_m)(1^{\ot 2} \ot x^{\wedge (\beta+[i])}) 
\ot {(x^*)}^{\wedge (\beta+[i])}\\
&= \sum_{i=1}^N \delta_{\beta_i,0}(\Ph^{-1}_{m-1}(b \ot {(x^*)}^{\wedge \beta}))
\left( (-1)^{\sum_{s=1}^i \beta_s} \left[
\left( \prod_{s=1}^i q_{s,i}^{\beta_s} \right) x_i \ot 1 -
\left( \prod_{s=i}^N q_{i,s}^{\beta_s} \right) \ot x_i \right] \ot x^{\wedge \beta} \right)\\
& \hspace{5.3in} \ot {(x^*)}^{\wedge (\beta+[i])}.
\end{split}
\end{equation*}
Using the $A^e$-linearity of $\Ph^{-1}_{m-1}(b \ot (x^*)^{\wedge \beta})$ 
and the definition of $\Ph$, the above expression may be rewritten to show
that $d_m^*(b\ot (x^*)^{\wedge\beta})$ is equal to 
\begin{equation}
\label{formula for d_m^* for B}
\sum_{i=1}^N \delta_{\beta_i,0} (-1)^{\sum_{s=1}^i \beta_s}
\left[
 \left( \prod_{s=1}^i q_{s,i}^{\beta_s} \right) x_i  b -
 \left( \prod_{s=i}^N q_{i,s}^{\beta_s} \right) b x_i  \right] 
\ot {(x^*)}^{\wedge (\beta+[i])}.
\end{equation}
We will use these expressions for the differentials in the sequel.

Let ${\bf q}$ be a tuple of scalars as in the introduction.
We define the {\em quantum exterior algebra} 
\begin{equation}\label{qea}
  \Wedge _{\bf q}(V) = T(V)/(x_ix_j+q_{i,j}x_jx_i \mid 1\leq i,j\leq N),
\end{equation}
where $T(V)$ is the tensor algebra on $V$.
Note that the relations corresponding to $i=j$ are
$2x_i ^2 =0$ since $q_{ii}=1$, that is  $x_i^2=0$ in $\Wedge_{\bf q}(V)$.
It follows that the dimension of $\Wedge _{\bf q}(V)$ as a (graded) vector space is
the same as that of the exterior algebra $\Wedge(V)$.

Let ${\bf q}^{-1}$ denote the tuple consisting of all inverses of components of ${\bf q}$,
that is, ${\bf q}^{-1} = \{ q_{i,j}^{-1}\}_{1\leq i,j\leq N}$.
Using this notation, we have the following cup product theorem.
Our main applications are to the cases $B=S_{\bf q}(V)$
and $B=S_{\bf q}(V)\rtimes G$ for some finite group $G$ in the next section.

\begin{theorem}
\label{cup product HH(A,B)}
Let $B$ be an $S_{\bf q}(V)$-bimodule algebra. Then
the Hochschild cohomology algebra $\HHD (S_{\bf q}(V),B)$ is a subquotient algebra of
the tensor product $B \ot \Wedge _{\bf q^{-1}}(V^*)$, that is, the cup product
on $\HHD (S_{\bf q}(V),B)$ descends from the product on $B \ot \Wedge _{\bf q^{-1}}(V^*)$.
\end{theorem}

\begin{proof}
Letting $A=S_{\bf q}(V)$,
by Remark~\ref{complement}, there is a chain map
$\psi: A^{\ot (m+2)}\rightarrow A^e\ot \Wedge^m (V)$ such that $\psi\phi$ is the
identity map.
(We will not need an explicit formula for $\psi$.) 

We will use the cup product as defined on the bar complex,
together with the chain maps $\phi$ and $\psi$: For two cocycles
$\mu$ and $\nu$ defined on $A^e\ot \Wedge (V)$, in degrees $m$ and
$n$ respectively, their cup product is given by
$$
   \phi^*(\psi^*(\mu)\smile \psi^*(\nu)).
$$
This function is defined  by its action on all elements of the form  
 $1^{\ot 2}\ot x_{j_1}\wedge\cdots\wedge x_{j_{m+n}}$, which we calculate next.
Let $S_{m,n}$ denote the set of all $m,n$-shuffles,
that is all $\rho\in S_{m+n}$ for which
$\rho(1)<\cdots <\rho(m)$ and $\rho(m+1)<\cdots <\rho(m+n)$.
Note that these shuffles form a set of coset representatives of
the subgroup $S_m\times S_n$ of $S_{m+n}$, where $S_m$ acts on 
$\{1,\ldots, m\}$ and $S_n$ acts on $\{m+1,\ldots,m+n\}$.
Thus we may write each $\pi\in S_{m+n}$ as $\pi =\rho\sigma\tau$ where
$\rho\in S_{m,n}$, $\sigma\in S_m$, $\tau\in S_n$. 
By Lemma \ref{sigma-tau}, writing $\rho\sigma\tau = (\rho\sigma\rho^{-1})(\rho\tau\rho^{-1})
\rho$, we have
\begin{eqnarray*}
&& \hspace{-.8cm} \phi^*(\psi^*(\mu)\smile\psi^*(\nu))
    (1^{\ot 2}\ot x_{j_1}\wedge \cdots\wedge x_{j_{m+n}}) \\
&=& (\psi^*(\mu)\smile\psi^*(\nu))\left(\sum_{\pi\in S_{m+n}}
   (\sgn\pi) q_{\pi}^{j_1,\ldots,j_{m+n}} \ot x_{j_{\pi(1)}}\ot\cdots\ot
    x_{j_{\pi(m+n)}} \ot 1 \right)\\
&=& \sum_{\pi\in S_{m+n}} (\sgn \pi) q_{\pi}^{j_1,\ldots,j_{m+n}} \psi^*(\mu)
   (1 \ot x_{j_{\pi(1)}}\ot\cdots\ot x_{j_{\pi(m)}} \ot 1) \\
&& \hspace{3in} \cdot \ \psi^*(\nu) 
   (1 \ot x_{j_{\pi(m+1)}}\ot\cdots\ot x_{j_{\pi(m+n)}} \ot 1)\\
&=& \sum_{\rho\in S_{m,n}} \ \sum_{\sigma\in S_m, \tau\in S_n}
   (\sgn \rho\sigma\tau) q_{\rho\sigma\tau} ^{j_1,\ldots,j_{m+n}} \psi^*(\mu)
   (1 \ot x_{j_{\rho\sigma\tau(1)}}\ot\cdots\ot x_{j_{\rho\sigma\tau(m)}} \ot 1)\\ 
&& \hspace{3in} \cdot \ \psi^*(\nu)(1 \ot x_{j_{\rho\sigma\tau(m+1)}}\ot
  \cdots\ot x_{j_{\rho\sigma\tau(m+n)}} \ot 1)\\
&=&\sum_{\rho\in S_{m,n}} (\sgn\rho) q_{\rho}^{j_1,\ldots,j_{m+n}} \left(\sum_{\sigma\in S_m}
   (\sgn\sigma) q_{\rho\sigma\rho^{-1}}^{j_{\rho(1)},\ldots,j_{\rho(m)}} \psi^*(\mu)
    (1 \ot x_{j_{\rho\sigma(1)}}\ot\cdots\ot x_{j_{\rho\sigma(m)}} \ot 1)\right)\\
   && \hspace{2.8cm}\cdot \left(\sum_{\tau\in S_n} (\sgn\tau) q_{\rho\tau\rho^{-1}}^{j_{\rho(m+1)},
          \ldots,j_{\rho(m+n)}} \psi^*(\nu) (1 \ot x_{j_{\rho\tau(m+1)}}\ot\cdots\ot
           x_{j_{\rho\tau(m+n)}} \ot 1)\right)\\
&=& \sum_{\rho\in S_{m,n}}
    (\sgn\rho) q_{\rho}^{j_1,\ldots,j_{m+n}} 
   \phi^*\psi^*(\mu)(1^{\ot 2} \ot x_{j_{\rho(1)}}\wedge\cdots \wedge x_{j_{\rho(m)}})\\
&& \hspace{3in} \cdot \ \phi^*\psi^*(\nu)(1^{\ot 2} \ot x_{j_{\rho(m+1)}}\wedge \cdots\wedge x_{j_{\rho(m+n)}}).
\end{eqnarray*}
Since $\psi\phi$ is the identity map, 
$\phi^*\psi^*(\mu) = \mu$ and $\phi^*\psi^*(\nu)=\nu$.
Now replacing $\mu$ by $b\ot (x^*)^{\wedge\beta}$ and
$\nu$ by $b'\ot (x^*)^{\wedge\beta'}$, where $b,b'\in B$,
we see that only one
summand in the above can be nonzero, and that is the summand corresponding
to $\rho$ for which
$$
   x^{\wedge\beta} = x_{j_{\rho(1)}}\wedge\cdots\wedge x_{j_{\rho(m)}}
\ \ \ \mbox{ and } \ \ \
   x^{\wedge\beta'} = x_{j_{\rho(m+1)}}\wedge\cdots\wedge x_{j_{\rho(m+n)}}.
$$
For this summand, we obtain
$
  (\sgn\rho) q_{\rho}^{j_1,\cdots,j_{m+n}} bb'.
$
Comparing, we find also that
$$
  (bb' \ot (x^*)^{\wedge \beta}\wedge (x^*)^{\wedge\beta'})
   (1^{\ot 2}\ot x_{j_1}\wedge\cdots\wedge x_{j_{m+n}}) = 
 (\sgn\rho) q_{\rho}^{j_1,\cdots,j_{m+n}} bb',
$$
by permuting the exterior factors in $(x^*)^{\wedge \beta}\wedge (x^*)^{\wedge \beta'}$
via $\rho$ before applying the function, and by using the identity $q_{\rho^{-1}}
^{j_{\rho(1)},\ldots, j_{\rho(m+n)}} = (q_{\rho}^{j_1,\ldots,j_{m+n}})^{-1}$
(a consequence of Lemma \ref{sigma-tau}).
Therefore the product is as claimed.
\end{proof}
	
\begin{remark}
It is not necessary that the characteristic of $\k$ be 0 for 
Theorem~\ref{cup product HH(A,B)}: While Wambst's proof that complex
\eqref{label: free resolution} is exact requires characteristic 0,
a similar resolution may be obtained in positive characteristic by
using a twisted tensor product construction
\cite{BO} or general theory for Koszul algebras \cite{BG}. 
Our proof of Theorem~\ref{cup product HH(A,B)} uses only
resolution \eqref{label: free resolution} and its embedding into the bar resolution.
\end{remark}

\end{section}
\begin{section}
{The Hochschild cohomology algebra of $S_{\bf q}(V)\rtimes G$}
\label{gvsg}

Let $G$ be a finite group acting
on $A=S_{\bf q}(V)$ by graded algebra automorphisms.
We are interested in the cohomology of the skew group algebra $A\rtimes G$.
Since the characteristic of $\k$ is zero, this is known to be isomorphic
to the $G$-invariant subalgebra of the cohomology $\HHD(A, A\rtimes G)$ of
$A$ with coefficients in $A\rtimes G$. (See for example~\cite{S}.)
In this section we will compute this latter cohomology, 
$\HH^{\bu}(A, A \rtimes G) = \Ext_{A^e}^{\bu}(A, A \rtimes G)$,
in the case when $G$ acts diagonally on the basis $x_1, \ldots, x_N$ of $V$.
Note that each $g$-component $A_g$ is a (left) $A^e$-module (see Section~2).
Also note that since $A\rtimes G = \bigoplus_{g\in G} A_g$ as an $A^e$-module,
$$
   \Ext^{\bu}_{A^e}(A,A\rtimes G) \cong \bigoplus_{g\in G}
    \Ext^{\bu}_{A^e}(A,A_g).
$$
We will compute the summands $\Ext^{\bu}_{A^e}(A,A_g)$.

Fix $g \in G$.
In Section~\ref{B}, we applied the $\Hom$ functor
$\Hom_{A^e}(\cdot, B)$, for any $A$-bimodule $B$, 
to the $A^e$-resolution of $A$ in \eqref{label: free resolution},
and made appropriate identification to obtain the complex~\eqref{new Hom(resolution) for B}.
When we specialize this complex to $B=A_g$, we obtain
\begin{equation}
\label{new Hom(resolution) with G}
0 \xrightarrow{}  A_g \xrightarrow{d_1^*} 
A_g \ot \Wedge^1(V^*) \xrightarrow{d_2^*} 
A_g \ot \Wedge^2(V^*) \xrightarrow{} \cdots,
\end{equation}
where  formula~(\ref{formula for d_m^* for B}) yields that
$d_m^*((a\# g)\ot (x^*)^{\wedge\beta})$ is equal to 
\begin{equation}
\label{formula for d_m^* with G}
\sum_{i=1}^N
\delta_{\beta_i,0} (-1)^{\sum_{s=1}^i \beta_s} \left[
\left( \left( \prod_{s=1}^i q_{s,i}^{\beta_s} \right) x_ia - 
\left( \prod_{s=i}^N q_{i,s}^{\beta_s}\right) a(\lexp{g}{x_i}) \right) \#g \right]
\ot {(x^*)}^{\wedge(\beta+[i])},
\end{equation}
for all $a \in A$ and $\beta \in \O^N$ with $|\beta|=m-1$.

\begin{subsection}{Additive structure}

Suppose $G$ acts diagonally on the basis $x_1,\ldots,x_N$ of $V$, that is, there exist scalars
$\lambda_{g,i} \in \k$ such that 
$\lexp{g}{x_i} = \lambda_{g,i}x_i, \text{ for all } g \in G, i \in \{1,\ldots,N\}$.

For each $g \in G$, define
\begin{equation}\label{Cg}
C_g := \left\{ \gamma \in (\N \cup \{-1\})^N \mid 
\text{ for each } i \in \{1, \ldots, N\}, \;
\prod_{s=1}^N q_{i,s}^{\gamma_s} = \lambda_{g,i} \text{ or } \gamma_i = -1 \right\}.
\end{equation}
The following theorem gives the structure of the Hochschild cohomology as a 
graded vector space. Recall definition (\ref{qea}) of the quantum exterior algebra.

\begin{theorem}
\label{HH(A:G)}
If  $G$ acts diagonally on the
chosen basis of $V$, then 
$\HH^{\bu}(S_{\bf q}(V),S_{\bf q}(V)_g)$ is the graded vector subspace of 
$S_{\bf q}(V)_g \ot \Wedge_{{\bf q}^{-1}}(V^*)$ given by:
$$
\HH^m(S_{\bf q}(V),S_{\bf q}(V)_g) \cong \bigoplus_{\substack{\beta \in \O^N \\ |\beta| = m}} 
\bigoplus_{\substack{\alpha \in \N^N \\ \alpha - \beta \in C_g}}
\span_\k\{(x^\alpha\#g) \ot {(x^*)}^{\wedge \beta}\},
$$
for all $m \in \N, g \in G$.
Therefore, 
$$
\HH^m(S_{\bf q}(V),S_{\bf q}(V) \rtimes G) 
\cong  \bigoplus_{g \in G}
\bigoplus_{\substack{\beta \in \O^N \\ |\beta| = m}} 
\bigoplus_{\substack{\alpha \in \N^N \\ \alpha - \beta \in C_g}}
\span_\k\{(x^\alpha\#g) \ot {(x^*)}^{\wedge \beta}\},
$$
for all $m \in \N$, and $\HH^m(S_{\bf q}(V)\rtimes G)$ is its $G$-invariant subspace.
\end{theorem}

\begin{example}
Fix $g \in G$. We wish to describe $\HHD_g:=\HHD(S_{\bf q}(V), S_{\bf q}(V)_g)$
when $N=2$. We will work out two special cases. When $q_{1,2}$ is not a root of unity
and $\lambda_{g,1}$, $\lambda_{g,2}$ are not both equal to 1, we have
\begin{eqnarray*}
\HH^0_g & \cong & \{0\}\\
\HH^1_g & \cong & \{0\}\\
\HH^2_g & \cong & \span_\k\{ (1 \# g) \ot x_1^* \wedge x_2^*\}.
\end{eqnarray*}
For the second case, assume 
$q_{1,2}$ is simultaneously a primitive $\ell$th root of unity, a $\ell_1$th root
of $\lambda_{g,1}$, and a $\ell_2$th root of $\lambda_{g,2}^{-1}$. Also assume 
$q_{1,2} \not = \lambda_{g,2}$ and $q_{1,2}^{-1} \not = \lambda_{g,1}$. Then we have
\begin{eqnarray*}
\HH^0_g & \cong & \span_\k \{ x_1^{\alpha_1}x_2^{\alpha_2} \# g \mid 
\alpha_1, \alpha_2 \in \N, \ \ell \text{ divides both } \alpha_1-\ell_2 \text{ and } \alpha_2-\ell_1 \} \\
\HH^1_g & \cong & \span_\k \{ (x_1^{\alpha_1}x_2^{\alpha_2} \# g) \ot x_1^* \mid 
\alpha_1, \alpha_2 \in \N, \ \ell \text{ divides both } \alpha_2-\ell_1 \text{ and } \alpha_1-\ell_2-1 \}\\
&&\bigoplus \span_\k \{ (x_1^{\alpha_1}x_2^{\alpha_2} \# g) \ot x_2^* \mid 
\alpha_1, \alpha_2 \in \N, \ \ell \text{ divides both } \alpha_1-\ell_2 \text{ and } \alpha_2-\ell_1-1 \} \\
\HH^2_g & \cong &  \span_\k\{ (1 \# g) \ot x_1^* \wedge x_2^*\}\\
&&\bigoplus \span_\k \{(x_1^{\alpha_1}x_2^{\alpha_2} \# g) \ot x_1^* \wedge x_2^* \mid 
\alpha_1, \alpha_2 \in \N, \ \ell \text{ divides } \alpha_1-\ell_2-1 \text{ and } \alpha_2-\ell_1-1\}.
\end{eqnarray*}
Thus $\HH^0_g$ is the free module over $Z(S_{\bf q}(V))$ generated by $x_1^{\ell_2}x_2^{\ell_1}\# g$,
$\HH^1_g$ is the free $Z(S_{\bf q}(V))$-module generated by $(x_1^{\ell_2 +1}x_2^{\ell_1 }\# g)\ot 
x_1^*$ and $(x_1^{\ell_2}x_2^{\ell_1+1}\# g)\ot x_2^*$, and $\HH^2_g$ is the direct sum of the
$\k$-linear span of $(1\# g)\ot x_1^*\wedge x_2^*$ and the free $Z(S_{\bf q}(V))$-module
generated by $(x_1^{\ell_2+1}x_2^{\ell_1+1}\# g)\ot x_1^*\wedge x_2^*$. 
\end{example}

Define
\begin{equation}\label{C}
C := \left\{ \gamma \in (\N \cup \{-1\})^N \mid 
\text{ for each } i \in \{1, \ldots, N\}, \;
\prod_{s=1}^N q_{i,s}^{\gamma_s} =1 \text{ or } \gamma_i = -1 \right\}.
\end{equation}
Taking $G$ to be the trivial group with one element, we immediately obtain the following:

\begin{corollary}
\label{HH(A)}
$\HH^{\bu}(S_{\bf q}(V))$ is the graded vector subspace of 
$S_{\bf q}(V)\ot \Wedge_{{\bf q}^{-1}}(V^*)$ given by:
$$
\HH^m(S_{\bf q}(V)) \cong \bigoplus_{\substack{\beta \in \O^N \\ |\beta| = m}} 
\bigoplus_{\substack{\alpha \in \N^N \\ \alpha - \beta \in C}}
\span_\k\{x^\alpha \ot {(x^*)}^{\wedge \beta}\},
$$
for all $m \in \N$.
\end{corollary}

\begin{example}\label{N=2}
When $N=2$, the description of 
$\HHD(S_{\bf q}(V))$ simplifies considerably. In this case, if $q_{1,2}$ is not
a root of unity, then
\begin{eqnarray*}
\HH^0(S_{\bf q}(V)) & \cong & \k\\
\HH^1(S_{\bf q}(V)) & \cong & \span_\k\{ x_1 \ot x_1^*, x_2 \ot x_2^*\}\\
\HH^2(S_{\bf q}(V)) & \cong & \span_\k\{ 1 \ot x_1^* \wedge x_2^*, x_1x_2 \ot x_1^* \ot x_2^*\},
\end{eqnarray*}
and if $q_{1,2}$ is a primitive $\ell$th root of unity, $\ell\geq 2$, then
\begin{eqnarray*}
\HH^0(S_{\bf q}(V)) & \cong & \span_\k \{ x_1^{\alpha_1}x_2^{\alpha_2} \mid 
\alpha_1, \alpha_2 \in \N, \ \ell \text{ divides both } \alpha_1 \text{ and } \alpha_2 \} \\
\HH^1(S_{\bf q}(V)) & \cong & \span_\k \{ x_1^{\alpha_1}x_2^{\alpha_2} \ot x_1^* \mid 
\alpha_1, \alpha_2 \in \N, \ \ell \text{ divides both } \alpha_1-1 \text{ and } \alpha_2 \}\\
&&\bigoplus \span_\k \{ x_1^{\alpha_1}x_2^{\alpha_2} \ot x_2^* \mid 
\alpha_1, \alpha_2 \in \N, \ \ell \text{ divides both } \alpha_1 \text{ and } \alpha_2-1 \} \\
\HH^2(S_{\bf q}(V)) & \cong &  \span_\k\{ 1 \ot x_1^* \wedge x_2^*\}\\
&&\bigoplus \span_\k \{x_1^{\alpha_1}x_2^{\alpha_2} \ot x_1^* \wedge x_2^* \mid 
\alpha_1, \alpha_2 \in \N, \ \ell \text{ divides } \alpha_1-1 \text{ and } \alpha_2-1\}.
\end{eqnarray*}
Note that in the first case, the center of $S_{\bf q}(V)$ is $Z(S_{\bf q}(V))=\k$,
while in the second case, $Z(S_{\bf q}(V))$ is generated by $x_1^{\ell}$ and
$x_2^{\ell}$. Thus in either case, $\HH^0(S_{\bf q}(V))=Z(S_{\bf q}(V))$ as expected,
$\HH^1(S_{\bf q}(V))$ is the free $Z(S_{\bf q}(V))$-module generated by $x_1\ot x_1^*$ and
$x_2\ot x_2^*$, and $\HH^2(S_{\bf q}(V))$ is the direct sum of the $\k$-linear span
of $1\ot x_1^*\wedge x_2^*$ and the free $Z(S_{\bf q}(V))$-module generated by
$x_1x_2\ot x_1^*\wedge x_2^*$. Similar expressions may be obtained when $N\geq 3$,
but there are more cases to consider due to the many more parameters involved.
\end{example}

We introduce some notation and lemmas before proving Theorem~\ref{HH(A:G)}.

Fix $g\in G$.
For any $\alpha \in \N^N,
\beta \in \O^N$, and $i \in \{1, \ldots, N\}$, define
$$
\Omega_g(\alpha, \beta, i) := 
\begin{cases}
0, \text{ if } \displaystyle \prod_{s=1}^N q_{i,s}^{\alpha_s - \beta_s} = \lambda_{g,i}, \\
0, \text{ if } \beta_i = 1,\\
\varepsilon(\beta,i) \displaystyle \left( \prod_{s=1}^i q_{i,s}^{\alpha_s-\beta_s} - 
\lambda_{g,i} \prod_{s=i}^N q_{s,i}^{\alpha_s-\beta_s} \right), \text{ otherwise},
\end{cases}
$$
where $\varepsilon(\beta,i) = (-1)^{\sum_{s=1}^i \beta_s}$. 

Then, using formula \eqref{formula for d_m^* with G} for $d_m^*$ we see that
\begin{equation}
\label{new formula for d_m^* with G}
d_m^*((x^\alpha \#g) \ot {(x^*)}^{\wedge \beta}) =
\sum_{i=1}^N \Omega_g(\alpha, \beta, i) (x^{\alpha+[i]} \#g) \ot {(x^*)}^{\wedge(\beta+[i])},
\end{equation}
for all $\alpha \in \N^N$ and $\beta \in \O^N$ with
$|\beta|=m-1$.

For any $\gamma \in (\M)^N$ and $m \in \N$, define
$$
K_{g,\gamma}^m := \span_\k\{(x^\alpha \#g) \ot {(x^*)}^{\wedge \beta} \mid
\alpha \in \N^N, \ \beta \in \O^N, \ |\beta|=m, \text{ and } 
\alpha-\beta=\gamma\}.
$$
Let $K_{g,\gamma}^{\bu}$ denote the subcomplex of \eqref{new Hom(resolution) with G}
whose $m$th term is given by $K_{g,\gamma}^m$. That $K_{g,\gamma}^{\bu}$ is
indeed a subcomplex of \eqref{new Hom(resolution) with G} follows from
\eqref{new formula for d_m^* with G}. We immediately obtain:

\begin{lemma}
\label{grading with G}
The complex \eqref{new Hom(resolution) with G} admits a grading by $(\M)^N$. Precisely,
the $m$th term of complex \eqref{new Hom(resolution) with G} decomposes as
$\displaystyle \bigoplus_{\gamma \in (\M)^N} K_{g,\gamma}^m$, for all $m \in \N$.
\end{lemma}

We now separately handle the cases where $\gamma$ is or is not in the set $C_g$ defined in \eqref{Cg}.
For the proof of the next lemma, we will need some notation:
For any $\gamma \in (\M)^N$, define
%
$$
\|\gamma\|_g :=  \# \left\{ i \in \{1,\ldots,N\}
   \mid \prod_{s=1}^N q_{i,s}^{\gamma_s} \neq \lambda_{g,i}
\text{ and } \gamma_i \neq -1\right\} .
$$

\begin{lemma}
\label{acyclic with G}
Let $\gamma \in ((\M)^N \backslash C_g)$. Then,
the subcomplex $K_{g,\gamma}^{\bu}$ of \eqref{new Hom(resolution) with G} is acyclic.
\end{lemma}
\begin{proof}
We will show that the identity chain map of the 
complex $K_{g,\gamma}^{\bu}$ is nullhomotopic. First,
define the following scalar. Let $\alpha \in \N^N, \beta \in \O^N$,
and $i \in \{1, \ldots N\}$. Define
$$
\omega_g(\alpha, \beta, i) := 
\begin{cases}
0, \text{ if } \displaystyle \prod_{s=1}^N q_{i,s}^{\alpha_s - \beta_s} = \lambda_{g,i}, \\
0, \text{ if } \alpha_i = 0,\\
0, \text{ if } \beta_i = 0,\\
\Omega_g(\alpha-[i], \beta-[i], i)^{-1}, \text{ otherwise}.
\end{cases}
$$
Now, fix $m \in \N$ and suppose that $|\beta|=m$ and
$\alpha-\beta = \gamma$. Define
$$
h_m : K_{g,\gamma}^m \to K_{g,\gamma}^{m-1}
$$
by
$$
h_m((x^\alpha \#g) \ot {(x^*)}^{\wedge \beta}) := \frac{1}{\|\gamma\|_g}
\sum_{i=1}^N \omega_g(\alpha,\beta,i) (x^{\alpha-[i]} \#g) \ot {(x^*)}^{\wedge(\beta-[i])}.
$$
Note that $\|\gamma\|_g \neq 0$, as $\gamma \not \in C_g$.

We contend that
$$
(h_{m+1} \circ d_{m+1}^*
+ d_m^* \circ h_m) ((x^\alpha \# g) \ot {(x^*)}^{\wedge\beta}) = (x^\alpha \# g) \ot {(x^*)}^{\wedge\beta}.
$$

Note that the lemma is proved if this equality is sustained.
The left hand side of this equality is equal to
\begin{equation*}
\begin{split}
&\frac{1}{\|\gamma\|_g} 
\sum_{i=1}^N \sum_{j=1}^N
\Omega_g(\alpha,\beta,i) \omega_g(\alpha+[i],\beta+[i],j)
(x^{\alpha+[i]-[j]}\#g) \ot {(x^*)}^{\wedge(\beta+[i]-[j])}\\
& \qquad \qquad \qquad +
\frac{1}{\|\gamma\|_g} 
\sum_{j=1}^N \sum_{i=1}^N
\omega_g(\alpha,\beta,j) \Omega_g(\alpha-[j],\beta-[j],i)
(x^{\alpha+[i]-[j]}\#g) \ot {(x^*)}^{\wedge(\beta+[i]-[j])}\\
=& \frac{1}{\|\gamma\|_g}
\sum_{i=1}^N \big[ \Omega_g(\alpha,\beta,i) \omega_g(\alpha+[i],\beta+[i],i)
+ \omega_g(\alpha,\beta,i) \Omega_g(\alpha-[i],\beta-[i],i) \big]
(x^{\alpha}\#g) \ot {(x^*)}^{\wedge \beta}\\
&\qquad \qquad \qquad +
\frac{1}{\|\gamma\|_g}
\sum_{i\neq j} \big[ \Omega_g(\alpha,\beta,i) \omega_g(\alpha+[i],\beta+[i],j)\\ 
&\qquad \qquad \qquad \qquad \qquad \qquad + \omega_g(\alpha,\beta,j) \Omega_g(\alpha-[j],\beta-[j],i) \big]
(x^{\alpha+[i]-[j]}\#g) \ot {(x^*)}^{\wedge (\beta+[i]-[j])}\\
\end{split}
\end{equation*}

It follows from the definition of $\Omega_g$ and $\omega_g$ that
$$
\sum_{i=1}^N \Omega_g(\alpha,\beta,i) \omega_g(\alpha+[i],\beta+[i],i)
+\omega_g(\alpha,\beta,i) \Omega_g(\alpha-[i],\beta-[i],i) = \|\gamma\|_g.
$$
Therefore, it only remains to show that
$$
\Omega_g(\alpha,\beta,i) \omega_g(\alpha+[i],\beta+[i],j) +
\omega_g(\alpha,\beta,j) \Omega_g(\alpha-[j],\beta-[j],i) = 0,
$$
whenever $i \neq j$. To this end, define
$$
\Xi_g(i,j) := \Omega_g(\alpha,\beta,i) \omega_g(\alpha+[i],\beta+[i],j)
$$
and
$$
\Xi'_g(i,j) := \omega_g(\alpha,\beta,j) \Omega_g(\alpha-[j],\beta-[j],i).
$$

Suppose $i \neq j$. Then
it is clear from definition of $\Omega_g$ and $\omega_g$ that
$\Xi_g(i,j)$ and $\Xi'_g(i,j)$ are simultaneously zero or nonzero, so
we may assume that they are nonzero. In this case, we have
$$
\Xi_g(i,j) = \varepsilon(\beta,i) \varepsilon(\beta+[i]-[j],j)
\left( \prod_{s=1}^i q_{i,s}^{\alpha_s-\beta_s} - 
\lambda_{g,i} \prod_{s=i}^N q_{s,i}^{\alpha_s-\beta_s} \right)
\left( \prod_{s=1}^j q_{j,s}^{\alpha_s-\beta_s} - 
\lambda_{g,j} \prod_{s=j}^N q_{s,j}^{\alpha_s-\beta_s} \right)^{-1}
$$
and
$$
\Xi'_g(i,j) = \varepsilon(\beta-[j],i) \varepsilon(\beta-[j],j)
\left( \prod_{s=1}^j q_{j,s}^{\alpha_s-\beta_s} - 
\lambda_{g,j} \prod_{s=j}^N q_{s,j}^{\alpha_s-\beta_s} \right)^{-1}
\left( \prod_{s=1}^i q_{i,s}^{\alpha_s-\beta_s} - 
\lambda_{g,i} \prod_{s=i}^N q_{s,i}^{\alpha_s-\beta_s} \right).
$$

Therefore, the desired equality $\Xi_g(i,j) = - \Xi'_g(i,j)$ is equivalent to the
equality
$$
\varepsilon(\beta-[j],i) \varepsilon(\beta-[j],j)
=- \varepsilon(\beta,i) \varepsilon(\beta+[i]-[j],j),
$$
which may be easily verified using the definition of $\varepsilon$,
the corresponding conditions under which $\Omega_g$ and $\omega_g$
are nonzero, and the condition $i \neq j$.
\end{proof}

\begin{proof}[Proof of Theorem \ref{HH(A:G)}]
Observe that the restriction of $d^*_{\bu}$ to $K_{g,\gamma}^{\bu}$
is zero, for all $\gamma \in C_g$. The theorem now 
follows immediately from Lemmas \ref{grading with G} and \ref{acyclic with G}. 
\end{proof}

\end{subsection}

\begin{subsection}{Cup product}

Assume that $G$ acts diagonally
on the basis $x_1,\ldots,x_N$ of $V$. 
The following is immediate from Theorem~\ref{cup product HH(A,B)} when we 
put $B=S_{\bf q}(V)\rtimes G$.
Recall the definition of $C_g$ in (\ref{Cg}) and of $\Wedge_{{\bf q}^{-1}}(V^*)$
given by (\ref{qea}).

\begin{theorem}
\label{cup product with G}
$\HHD (S_{\bf q}(V), S_{\bf q}(V) \rtimes G)$ is a subquotient algebra of 
$(S_{\bf q}(V) \rtimes G) \ot \Wedge_{{\bf q}^{-1}} (V^*)$.
Thus 
the cup product on $\HHD(S_{\bf q}(V), S_{\bf q}(V) \rtimes G)$ is given by
$$
 ( (x^{\alpha} \# g)\ot (x^*)^{\wedge \beta} ) \smile
  ((x^{\alpha '} \# h) \ot (x^*)^{\wedge \beta'})
=
    (x^{\alpha} \# g)(x^{\alpha '} \# h) \ot (x^*)^{\wedge\beta}\wedge (x^*)^{\wedge\beta '}
$$
for all $g,h\in G$, $\alpha,\alpha '\in \N^N$ and $\beta,\beta'\in \{0,1\}^N$ 
such that $\alpha-\beta\in C_g$ and $\alpha'-\beta'\in C_h$.
Moreover, $\HHD(S_{\bf q}(V)\rtimes G)$ is the $G$-invariant subalgebra of
$\HHD(S_{\bf q}(V),S_{\bf q}(V)\rtimes G)$.
\end{theorem}

\begin{remark}
Note that the above product is zero when the 
supports of $\beta$ and $\beta'$ intersect nontrivially. Furthermore, we 
understand the above product to be zero (i.e., a coboundary) if
$\alpha + \alpha ' -(\beta +\beta ')$ is not in $C_{gh}$.
\end{remark}

Recall the definition of $C$ in \eqref{C}.
Taking $G$ to be the trivial group with one element, we immediately obtain the following:

\begin{corollary}
$\HHD (S_{\bf q}(V))$ is a subquotient algebra of 
$S_{\bf q}(V)\ot \Wedge_{{\bf q}^{-1}} (V^*)$.
Thus 
the cup product on $\HHD(S_{\bf q}(V))$ is given by
$$
 ( x^{\alpha}\ot (x^*)^{\wedge \beta} ) \smile
  (x^{\alpha '}\ot (x^*)^{\wedge \beta'})
=
    x^{\alpha}x^{\alpha '} \ot (x^*)^{\wedge\beta}\wedge (x^*)^{\wedge\beta '}
$$
for all $\alpha,\alpha '\in \N^N$ and $\beta,\beta'\in \{0,1\}^N$ 
such that $\alpha-\beta\in C$ and $\alpha'-\beta'\in C$.
\end{corollary}

\begin{remark}
As before, the above product is zero when the 
supports of $\beta$ and $\beta'$ intersect nontrivially. Furthermore, we 
understand the above product to be zero (i.e., a coboundary) if
$\alpha + \alpha ' -(\beta +\beta ')$ is not in $C$. 
\end{remark}

\end{subsection}

\end{section}




\begin{thebibliography}{KKZ}

\bibitem[AS]{AS} N.\ Andruskiewitsch and H.-J.\ Schneider,
\textit{Pointed Hopf algebras}, in: New directions in Hopf algebras,
MSRI Publ., \textbf{43}, Cambridge Univ.\ Press, Cambridge, 2002.

\bibitem[A]{A} R.\ Anno, 
\textit{Multiplicative structure on the Hochschild cohomology of 
crossed product algebras},
\texttt{arXiv:math.QA/0511396}.

\bibitem[BB]{BB} Y.\ Bazlov and A.\ Berenstein, 
\textit{Noncommutative Dunkl operators and braided Cherednik algebras},
Sel. Math. New Ser. \textbf{14} (2009), no.\ 3--4, 325--372.

\bibitem[BGS]{BGS} A.\ Beilinson, V.\ Ginzburg and V.\ Soergel,
\textit{Koszul duality patterns in representation theory},
J.\ Am.\ Math.\ Soc.\ \textbf{9} (1996), no.\ 2, 473--527.

\bibitem[BO]{BO} P.\ A.\ Bergh and S.\ Oppermann, 
\textit{Cohomology of twisted tensor products}, 
J.\ Algebra \textbf{320} (2008), no. 8, 3327--3338.

\bibitem[BG]{BG} A.\ Braverman and D.\ Gaitsgory, \textit{Poincar\'e-Birkhoff-Witt Theorem
for quadratic algebras of Koszul type}, J.\ Algebra \textbf{181} (1996), no.\ 2, 315--328.

\bibitem[F]{F} M.\ Farinati,
\textit{Hochschild duality, localization, and smash products},
J.\ Algebra \textbf{284} (2005), no.\ 1, 415--434.

\bibitem[G]{G} M.\ Gerstenhaber, 
\textit{The cohomology structure of an associative ring},
Ann.\ of Math.\ \textbf{78} (1963), no.\ 2, 267--288.

\bibitem[GG]{GG} J.\ A.\ Guccione and J.\ J.\ Guccione, 
\textit{Hochschild and cyclic homology of Ore extensions and some
examples of quantum algebras}, 
K-Theory \textbf{12} (3) (1997), 259--276.

\bibitem[Gu]{Gu} A.\ Guichardet, 
\textit{Homologie de Hochschild des d\'eformations quadratiques d'alg\'ebres
de polyn\^omes},
Comm.\ Algebra \textbf{26} (12) (1998), 4309--4330.

\bibitem[GK]{GK} V.\ Ginzburg and D.\ Kaledin,
\textit{Poisson deformations of symplectic quotient singularities},
Adv.\ Math.\ \textbf{286} (2004), no.\ 1, 1--57.


\bibitem[KKZ]{KKZ} E.\ Kirkman, J.\ Kuzmanovich, and J.\ J.\ Zhang,
\textit{Shephard-Todd-Chevalley Theorem for skew polynomial rings},
Alg. Rep. Theory \textbf{13} (2010), no. 2, 127--158. 

\bibitem[M]{M} Yu.\ I.\ Manin, \textit{Quantum groups and non-commutative
geometry}, CRM Universit\'e de Montr\'eal, 1988.

\bibitem[P]{P} S.\ Priddy, \textit{Koszul resolutions}, Trans.\ Amer.\ Math.\ Soc.\
\textbf{152} (1970), 39--60.

\bibitem[R]{R} L.\ Richard,
\textit{Hochschild homology and cohomology of some classical and
quantum noncommutative polynomial algebras},
J.\ Pure Appl.\ Algebra \textbf{187} (2004), 255--294. 

\bibitem[SW]{SW} A.\ V.\ Shepler and S.\ Witherspoon,
\textit{Hochschild cohomology and graded Hecke algebras},
Trans.\ Amer.\ Math.\ Soc.\ \textbf{360} (2008), no.\ 8, 3975--4005.

\bibitem[S]{S} D.\ \c{S}tefan, 
\textit{Hochschild cohomology on Hopf Galois extensions},
J.\ Pure Appl.\ Algebra \textbf{103} (1995), 221-233.

\bibitem[W]{W} M.\ Wambst, 
\textit{Complexes de Koszul quantiques},
Ann.\ Fourier, \textbf{43} (1993), no.\ 4, 1089--1156.


\end{thebibliography}
\end{document}